\documentclass[a4paper,12pt]{article}
\usepackage[hideboxes]{boxedeps}

\usepackage{latexsym}
\newcommand{\pic}[2]{\BoxedEPSF{#1 scaled #2}}

\def\Firstbraid{\pic{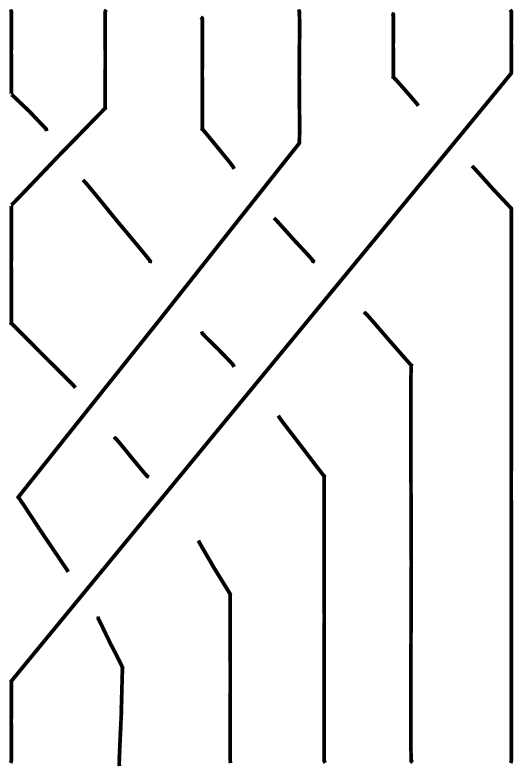} {350}}
\def\Secondbraid{\pic{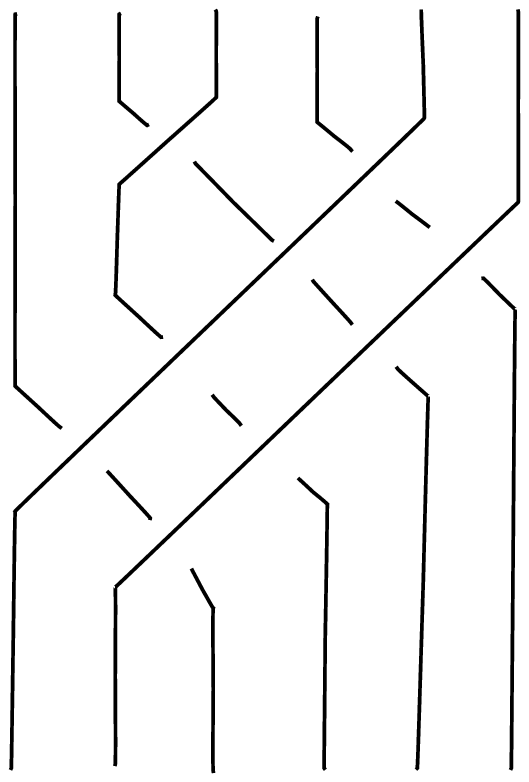} {350}}
\def\Firstsquare{\pic{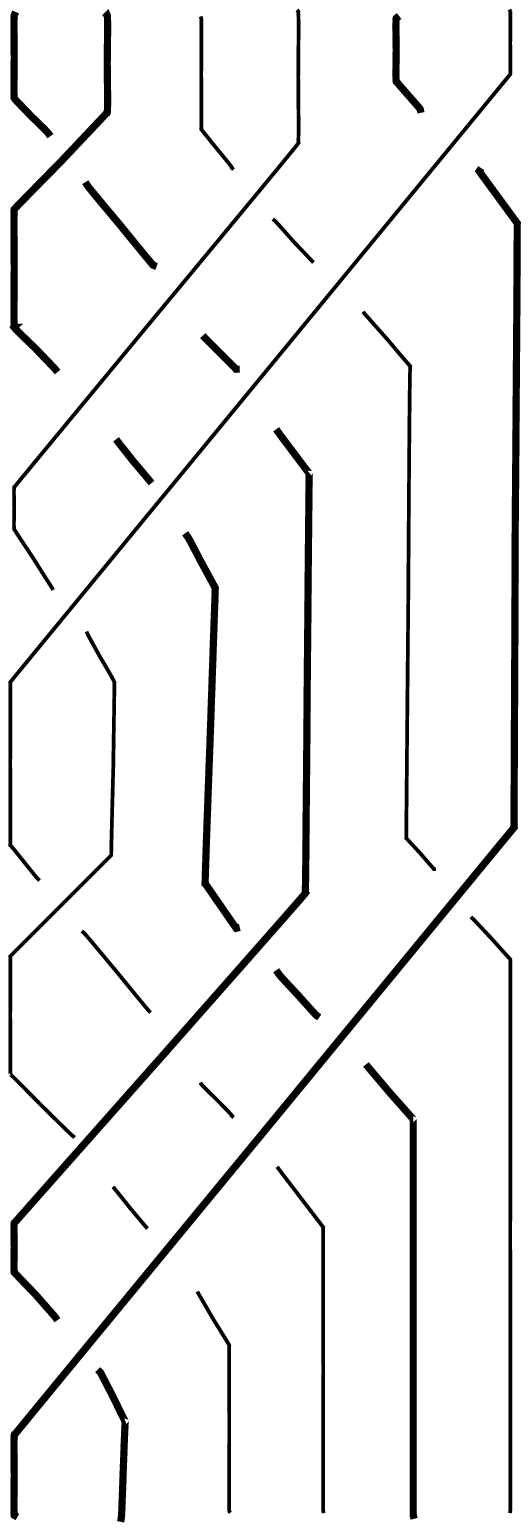} {300}}
\def\Secondsquare{\pic{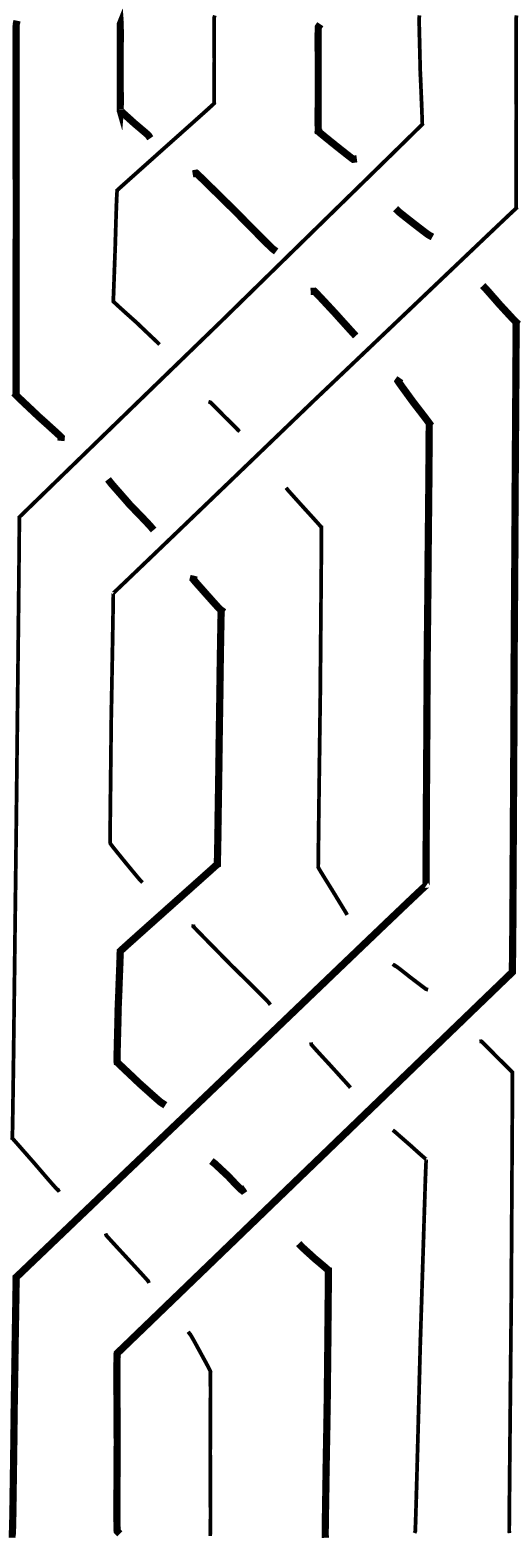} {300}}

\newcommand{\s}{\sigma}
\newcommand{\be}{\beta}

\newtheorem{theorem}{Theorem}

\newtheorem{lemma}[theorem]{Lemma}

\newenvironment{remark}{\par\smallskip%
\noindent\textbf{Remark.}\  }%
{\par\smallskip}

\newenvironment{notation}{\par\smallskip%
\noindent\textbf{Notation.}\  }%
{\par\smallskip}

\newenvironment{Def}{\par\smallskip%
\noindent\textbf{Definition.}\  }%
{\par\smallskip}

\newenvironment{proof}[1][{}]{\par\smallskip%
\noindent\textit{Proof #1: }\  }
{\hfill$\Box$\par\smallskip}

\begin{document}
\begin{center}
      {\Large\bf Conjugacy for positive permutation braids}

\bigskip
HUGH R. MORTON and RICHARD J. HADJI

\medskip
 {\it Department of Mathematical Sciences, University of Liverpool,\\
     Peach St, Liverpool, L69 7ZL, England. }\\
{\tt morton@liv.ac.uk\,;\,rhadji@liv.ac.uk}
\end{center}
\vskip 1cm


\begin{abstract}
Positive permutation braids on $n$ strings, which are defined to be  positive $n$-braids where each
pair of strings crosses at most once, form the elementary but non-trivial building blocks in many
studies of conjugacy in the braid groups. We consider conjugacy among these elementary braids which
close to knots, and show that those which close to the trivial knot or to the trefoil are all
conjugate. All such
$n$-braids with the maximum possible crossing number are also shown to be conjugate.

We note that conjugacy of these braids for $n\le5$ depends only on the crossing number. In
contrast, we exhibit two such braids on $6$ strings with $9$ crossings which are not conjugate but
whose closures are each isotopic to the $(2,5)$ torus knot.
\phantom{some stuff to separate the keywords from the abstract}

\end{abstract}

\section*{Introduction} The question of when two $n$-string braids are conjugate has aroused
interest over many years. Algorithms for comparing braids based on refinements of Garside's
algorithm,
\cite{EM,Meneses},
 can be used to settle this question in individual cases. A basic complexity measure in the
algorithms
 is the least number of  permutation braids needed to present a conjugate of
the given braid.  The simplest general case, when this number is 1,  reduces to deciding when two
positive permutation braids on $n$ strings are conjugate.

A necessary condition is that the corresponding permutations be conjugate, in other words the
permutations have the same cycle type.  In this investigation we shall restrict ourselves to the
case where the closure of the braid is a knot, and equivalently to those permutations in $S_n$ 
which are
$n$-cycles.

While any two such permutations are conjugate, the corresponding permutation braids need not be.
A sufficient condition is that the closures of the braids be isotopic {\em as closed braids}, in
other words the closed braids must be isotopic in the solid torus which is the complement of the
braid axis, \cite{Morton}.

We shall examine how far this condition follows from weaker necessary conditions on the braids.

\setcounter{theorem}{1}

\begin{theorem}  Positive permutation braids on $n$ strings which close to the unknot are
all conjugate.
\end{theorem}

\begin{theorem}  Positive permutation braids on $n$ strings which close to the trefoil
are all conjugate.
\end{theorem}

\begin{theorem}  Positive permutation braids on $n$ strings which close to the same knot
are all conjugate, when $n\le5$.
\end{theorem}

We also prove a general result in theorem \ref{large} about conjugacy of such braids which have the
largest possible number of crossings. On the other hand in theorem \ref{nonconj}  we  exhibit
two
$6$-string positive permutation braids which close to the (2,5) torus knot but are not conjugate.
These are constructed along the lines of Murasugi and Thomas' original example of non-conjugate
positive braids with isotopic closure,
\cite{Thomas}.

Further simple non-conjugacy results in theorem \ref{trefoil} give  a range of non-conjugate
positive braids closing to the trefoil, in  contrast to theorem \ref{2}.

Some of our results were first noted in \cite{Hadji} by the second author. There has also been a
recent exploration by Elrifai and Benkhalifa \cite{EB} for small values of $n$ without restrictions
on the cycle type of the permutation.

The techniques used in this paper to prove non-conjugacy are very direct; more subtle techniques,
such as Fiedler's Gauss sum invariants \cite{Fiedler}, may be used in more difficult cases, or
applications of the algorithm of Franco and Meneses \cite{Meneses}. Hall has examples coming from
the realms of dynamical systems of positive permutation braids on 12 or more strings which
are believed not to be conjugate to their reverse \cite{Hall}. In such cases none of the techniques
used here can be applied to establish non-conjugacy. 
\setcounter{theorem}{0}

\section{Permutation braids}
We shall use Artin's classical description of the group $B_n$ of braids on $n$ strings in terms of
elementary generators $\s_i$  for $i=1,\ldots,n-1$ with the relations:
\begin{enumerate}
  \item $\s_i\s_j=\s_j\s_i$ for $|i-j|\ge 2$,
  \item $\s_i\s_{i+1}\s_i=\s_{i+1}\s_i\s_{i+1}$ for $1\le i\le n-2$.
\end{enumerate}

There are two simple homomorphisms from the group $B_n$ which give initial constraints on conjugacy.

\begin{itemize}
\item
The  homomorphism $\varphi:B_n\to S_n$ defined on the generators by $\varphi(\sigma_i)=(i\;i+1)$
determines a permutation $\pi=\varphi(\beta)$ in which $\pi(j)$ gives the endpoint of the string of
$\beta$ which begins at $j$.

\item
The homomorphism $wr:B_n\to {\bf Z}$ defined  by $wr(\sigma_i)=1$,  counts the
{\em writhe} or `algebraic crossing number' of a braid. 
\end{itemize}

Two conjugate braids in $B_n$ must then
have the same writhe, since $\bf Z$ is abelian, as well as having permutations of the same cycle
type.

\begin{Def}
A {\em positive braid} is an element of $B_n$ which can be
written as a word in positive powers of the generators $\{\s_i\}$,
without use of the inverse elements $\s_i^{-1}$.
\end{Def}
For positive braids, the writhe is simply the number of crossings in the braid.

\begin{Def}
A braid $\beta$ is called a {\em positive permutation braid} if
it is a positive braid such that  no pair of strings
cross more than once.
\end{Def}

\begin{notation}
We denote the set of positive braids and positive permutation
braids in $B_n$ by $B_n^+$ and $S_n^+$ respectively.
\end{notation}

This definition of positive permutation  braids was first used by Elrifai in
\cite{Elrifai, EM}, where they were shown to correspond exactly to permutations. Explicitly, the
 homomorphism $\varphi$ restricts to a bijection from the set $S_n^+$ of  positive permutation
braids to $S_n$. They were also identified by Elrifai with the set of  {\em initial segments} of
Garside's fundamental braid $\Delta_n$.

It should be noted  that the explicit braid word for a positive permutation braid is generally
not unique. For example,  the permutation $(1423)$ can be
represented  in $S_4^+$ by braid words $\s_1\s_2\s_3\s_1\s_2$ and
$\s_2\s_1\s_3\s_2\s_3$.
Consequently some authors choose to label  permutation
braids simply by the corresponding permutation in $S_n$.

The number of components of the closure $\hat{\beta}$ of a braid $\beta$, constructed by
identifying the initial points with the end points, is the number of cycles in the cycle type of
the permutation
$\varphi(\beta)$.
 In this paper we restrict attention to braids which close to knots, and hence we
shall only look at the $(n-1)!$ permutation braids whose permutation is a single
$n$-cycle.

\section{Conjugacy results}

Suppose that two braids $\beta$ and $\gamma$ are conjugate in $B_n$. Then their closures are
isotopic as links in the complement of the braid axis, and so they are  certainly isotopic  in
$S^3$.

Where $\beta$ and $\gamma$ are positive  they must have the same number of crossings, because they
have the same writhe. Even if they are not conjugate, two positive braids in $B_n$ which close to
isotopic knots must have the same number of crossings. 

\begin{lemma}
If $\be,\gamma\in B_n^+$ and $\hat{\be},\hat{\gamma}$ are isotopic knots then $wr(\be)=wr(\gamma)$.
\end{lemma}
\begin{proof}
Suppose that the knot $\hat{\be}$ has genus $g$.
 The closure of a  non-split positive braid $\beta\in B_n$ is always a fibred link or
knot. The surface found from $\hat{\be}$ by Seifert's algorithm is a fibre surface, and has minimal
genus
$g$. Its Euler characteristic
$\chi=1-2g$ satisfies
$1-\chi=c-(n-1)$ where
$c=wr(\be)$ is the number of crossings in $\be$. Hence  
$wr(\be)=(n-1)+2g=wr(\gamma)$, since $\hat{\gamma}$ is isotopic to $\hat{\be}$ and so also has
genus $g$.
\end{proof}

Consequently if a positive $n$-braid closes to the unknot then it must have exactly $n-1$
crossings.  If it  closes to the trefoil knot, which has genus 1, then it must have $n+1$ crossings.
We now show that positive permutation $n$-braids which close to either of these knots are determined
up to conjugacy.

Explicitly we have the following  results.
\begin{theorem} \label{1}

Any positive permutation $n$-braid $\beta$ which closes to the unknot is conjugate to
$\sigma_1\sigma_2
\cdots\sigma_{n-1}$.
\end{theorem}

\begin{theorem}\label{2}

Any positive permutation $n$-braid $\beta$ which closes to the trefoil is conjugate to
$\sigma_1^3\sigma_2\cdots\sigma_{n-1}$.
\end{theorem}

\begin{proof} [of theorem \ref{1}]

 Each generator $\sigma_i$ must appear at least
once in $\beta$, otherwise its closure is disconnected. Since its closure has genus $0$ the
braid $\beta$ has $n-1$ crossings, and so  each generator appears exactly once.

It is enough to manipulate the braid cyclically, as such manipulations can be realised as
conjugacies. We can represent $\beta$ up to conjugacy
by writing the generators $\s_1,\ldots,\s_{n-1}$ in the appropriate order around a circle. Each
generator appears exactly once. To prove the theorem we  use the braid  commutation relations to
rearrange the generators in ascending order round the circle.

 Assume by induction on
$j$ that the generators
$\s_1,\ldots,\s_j$ occur consecutively in order. Then any generator $\s_k$  lying on the circle
between $\s_j$ and
$\s_{j+1}$  has $k>j+1$. These  generators then commute with each of $\s_1,
\ldots,\s_j$, and can be moved past them to leave
$\s_{j+1}$ immediately after
$\s_j$.
The process finishes when all generators are in consecutive order.
\end{proof}
 
\begin{proof} [of theorem \ref{2}]

Again represent  generators on a circle.
Each generator $\s_1,\ldots,\s_{n-1}$ must occur at least once, otherwise the closed braid splits.
Since the trefoil has genus $1$ the braid has $n+1$ crossings. So either two generators $\s_i$ and
$\s_j$ each occur twice, or one,
$\s_i$ say, occurs 3 times, and the other generators occur once only.

Either $\s_{i-1}$ or $\s_{i+1}$ must occur between two occurrences of $\s_i$ in  $\beta$,
otherwise it can be rewritten with two
consecutive occurrences of $\s_i$. This is not possible for a permutation braid, since pairs of
strings cross at most once.
If $\s_i$
occurs three times then $\s_{i-1}$  lies between one pair of occurrences of $\s_i$  and $\s_{i+1}$
between the other pair. We can then move all generators except $\s_{i+1}$ past these last two
occurrences of $\s_i$ to write $\beta$ with a consecutive sequence $\s_i\s_{i+1}\s_i$. Change this to
$\s_{i+1}\s_i\s_{i+1}$ by the braid relation to write $\beta$ with $\s_i$ and $\s_{i+1}$ each
appearing twice.

We may thus assume that two generators $\s_i$ and $\s_j$ each occur twice in $\beta$, with $j>i$. If
$j>i+1$ then
$\s_{i+1}$ occurs only once. We can then collect all generators $\s_k$ with $k>i+1$ at the
two ends of the braid word, and combine them at the end of the word by cycling so as to write a
conjugate braid in the form $AB$ where  $B$ is a product of generators $\s_k $ with
$k>i+1$, and includes
$\s_j$ twice, while  $A$ is a product with $k\le i+1$,  and includes $\s_i$ twice. The
closure of the braid then has three components and not one.

Hence $\beta$ must contain $\s_i$ and $\s_{i+1}$ twice each. Furthermore their occurrences must be
interleaved, otherwise we can cycle the braid and commute elements to separate it as a product of
generators $\s_k$ with $k\le i$ and those with $k>i$, and  its closure will again have three
components.

We shall prove, by induction on $i$, that any positive braid with two interleaved occurrences of
$\s_i$ and $\s_{i+1}$, and single occurrences of all other generators, is conjugate to
$\sigma_1^3\sigma_2\cdots\sigma_{n-1}$. 

 We can assume, by cycling, that the single occurrence of
$\s_{i-1}$ does not lie between the two occurrences of $\s_i$. We can move all further  generators
except
$\s_{i+1}$ past $\s_i$ so as to write
$\s_i\s_{i+1}\s_i$ consecutively.  The remaining occurrence  of $\s_{i+1}$ can be moved round the
circle past any other generator except the single $\s_{i+2}$. It can then be moved one way or other
round the circle to reach this block of three generators, giving either $\s_i\s_{i+1}\s_i\s_{i+1}$
or
$\s_{i+1}\s_i\s_{i+1}\s_i$. The braid relation then gives a consecutive block of either
$\s_{i}\s_i\s_{i+1}\s_i$ or $\s_i\s_{i+1}\s_i\s_{i}$.  The single $\s_{i-1}$ now lies between two
occurrences of
$\s_i$ on the circle. Any intervening generators commute with $\s_i$ and can be moved out  to leave
$\s_i\s_{i-1}\s_i$, which can be converted to $\s_{i-1}\s_i\s_{i-1}$. The cyclic braid now has two
interleaving occurrences of $\s_{i-1}$ and
$\s_i$.

 The result follows by induction on $i$, once we establish it for $i=1$.
In this case the argument above provides a block of either $\s_1\s_1\s_2\s_1$ or $\s_1\s_2\s_1\s_1$
on the circle. Since
$\s_1$ commutes with all generators except $\s_2$ we can move the right hand occurrences of $\s_1$
round the circle to give a block $\s_1\s_1\s_1\s_2$. The remaining generators can then be put in
ascending order as in the proof of theorem \ref{1}.
\end{proof}

A quick check on the possible values of the writhe for the $(n-1)!$ positive permutation braids
with $n\le4$ which close to a knot shows that in this range conjugacy is determined simply by
writhe, using theorems \ref{1} and \ref{2}. Positive permutation braids with $n+3$ crossings arise
first when $n=5$. A direct check on the corresponding braids shows that in this case too the
writhe is sufficient. 

\begin{theorem}\label{3}
   Positive permutation braids on $n$ strings which close to a knot
are  conjugate if and only if they have the same number of crossings, when $n\le5$.
\end{theorem}

Tables of these braids for $n=3,4,5$, and the corresponding permutations,  are included below.

When $n=3$ there are just two braids which both close to the unknot.
\begin{center}
\begin{tabular}{|c|c|c|}\hline
  Permutation & Braid word & Number of \\
              &            & crossings \\ \hline
  (123)       & $\s_2\s_1$ & 2 \\ \hline
  (132)       & $\s_1\s_2$ & 2 \\ \hline
\end{tabular}
\end{center}

When $n=4$ there
 are two conjugacy classes. The braids with writhe $3$ close to the unknot, and those with writhe
$5$ to the trefoil.
\begin{center}
\begin{tabular}{|c|c|c|}\hline
  Permutation & Braid word & Number of \\
              &            & crossings \\ \hline
  (1234)      & $\s_3\s_2\s_1$ & 3 \\ \hline
  (1243)      & $\s_2\s_1\s_3$ & 3 \\ \hline
  (1342)      & $\s_1\s_3\s_2$ & 3 \\ \hline
  (1432)      & $\s_1\s_2\s_3$ & 3 \\ \hline
  (1324)      & $\s_2\s_1\s_3\s_2\s_1$ & 5 \\ \hline
  (1423)      & $\s_1\s_2\s_1\s_3\s_2$ & 5 \\ \hline
\end{tabular}
\end{center}

When $n=5$ there are three conjugacy classes. The braids with writhe $4=n-1$  close to the unknot
and those with writhe $6$ close to the trefoil. Those with writhe
$8=n+3$  all close to the
$(2,5)$ torus knot.
\begin{center}
\begin{tabular}{|c|c|c|}\hline
  Permutation & Braid word                       & Number of \\
              &                                  & crossings \\ \hline
  (12345)     & $\s_4\s_3\s_2\s_1$                 & 4 \\ \hline
  (12354)     & $\s_3\s_2\s_1\s_4$                 & 4 \\ \hline
  (12453)     & $\s_2\s_1\s_4\s_3$                 & 4 \\ \hline
  (12543)     & $\s_2\s_1\s_3\s_4$                 & 4 \\ \hline
  (13452)     & $\s_1\s_4\s_3\s_2$                 & 4 \\ \hline
  (13542)     & $\s_1\s_3\s_2\s_4$                 & 4 \\ \hline
  (14532)     & $\s_1\s_2\s_4\s_3$                 & 4 \\ \hline
  (15432)     & $\s_1\s_2\s_3\s_4$                 & 4 \\ \hline
  (12435)     & $\s_3\s_2\s_4\s_3\s_2\s_1$         & 6 \\ \hline
  (12534)     & $\s_2\s_3\s_2\s_1\s_4\s_3$         & 6 \\ \hline
  (13245)     & $\s_2\s_1\s_4\s_3\s_2\s_1$         & 6 \\ \hline
  (13254)     & $\s_2\s_1\s_3\s_2\s_1\s_4$         & 6 \\ \hline
  (13524)     & $\s_3\s_2\s_1\s_4\s_3\s_2$         & 6 \\ \hline
  (14253)     & $\s_2\s_1\s_3\s_2\s_4\s_3$         & 6 \\ \hline
  (14352)     & $\s_1\s_3\s_2\s_4\s_3\s_2$         & 6 \\ \hline
  (14523)     & $\s_1\s_2\s_1\s_4\s_3\s_2$         & 6 \\ \hline
  (15342)     & $\s_1\s_2\s_3\s_2\s_4\s_3$         & 6 \\ \hline
  (15423)     & $\s_1\s_2\s_1\s_3\s_2\s_4$         & 6 \\ \hline
  (13425)     & $\s_2\s_3\s_2\s_1\s_4\s_3\s_2\s_1$ & 8 \\ \hline
  (14235)     & $\s_1\s_3\s_2\s_1\s_4\s_3\s_2\s_1$ & 8 \\ \hline
  (14325)     & $\s_2\s_1\s_3\s_2\s_4\s_3\s_2\s_1$ & 8 \\ \hline
  (15234)     & $\s_1\s_2\s_3\s_2\s_1\s_4\s_3\s_2$ & 8 \\ \hline
  (15243)     & $\s_1\s_2\s_1\s_3\s_2\s_4\s_3\s_2$ & 8 \\ \hline
  (15324)     & $\s_1\s_2\s_1\s_3\s_2\s_1\s_4\s_3$ & 8 \\ \hline
\end{tabular}
\end{center}

\bigskip
Having looked  among the closures of positive permutation
braids at knots with the smallest number of crossings, in theorems \ref{1} and \ref{2}, we now
turn briefly to those with the largest possible number.

The largest number of crossings in any positive permutation braid in $B_n$ is $\frac{1}{2}n(n-1)$,
which occurs for the fundamental half-twist braid $\Delta_n$. If the closure is to be a knot the
largest number of crossings is $\frac{1}{2}n(n-1)-[\frac{1}{2}(n-1)]$.

\begin{theorem}\label{large}
Every positive permutation braid with $\frac{1}{2}n(n-1)-[\frac{1}{2}(n-1)]$ crossings which closes
to a knot is conjugate to $\Delta_n \sigma_1^{-1}\sigma_2^{-1}\cdots\sigma_k^{-1}$ where
$k=[\frac{1}{2}(n-1)]$.
\end{theorem}

\begin{proof}
Take $n=2k+1$ or $n=2k+2$, so that $k=[\frac{1}{2}(n-1)]$, and let $\beta$ be a positive
permutation braid with $\frac{1}{2}n(n-1)-k$ crossings which closes
to a knot. Then $\beta$ has a complementary positive permutation braid $\gamma$ in $\Delta_n$, with
$\be\gamma=\Delta_n$. The braid $\gamma$ has $k$ crossings. Since $\beta=\Delta_n\gamma^{-1}$
closes to a knot the $k$ crossings in $\gamma^{-1}$ must be used to connect up the $k+1$ components
in $\hat\Delta_n$. Hence the $k$ generators  in $\gamma^{-1}$ must all be different. When $n=2k+2$
 the generator $\s_{k+1}$ cannot occur, since this connects two strings which are already in
the same component of $\hat\Delta_n$, and more generally  $\s_j$ and
$\s_{n-j}$ cannot both occur, for any $j$, as they both connect the same two components. In
particular  the generators $\s_k$
and
$\s_{k+1}$ cannot both occur when $n=2k+1$.

The generators in $\gamma$ then belong to two mutually commuting sets, those from $\s_1$ up to
$\s_k$ and those from $\s_{k+1}$ up to $\s_{n-1}$.  Write $\be$ on a circle, with one block of
generators together as $\Delta_n$ and then the $k$ generators of $\gamma^{-1}$. Move all
generators $ \s_j$ with $j>k$ to the extreme right in $\gamma^{-1}$ and then round the circle to the
left of $\Delta_n$. Now move them past $\Delta_n$, when each $\s_j$ is converted to $\s_{n-j}$.
Since $\s_j$ and $\s_{n-j}$ did not both occur in $\gamma$ we get a braid $\Delta_n\alpha^{-1}$
conjugate to $\be$ in which $\s_1,\ldots,\s_k$ each occurs exactly once in $\alpha$.

Following the method of theorem \ref{1} we can arrange the $k$ generators in $\alpha$ in any order
up to conjugacy, once we know how to move any generator $\s_j$ from the left to the right of
$\alpha^{-1}$ by conjugacy.  This can be done by taking it {\em twice} round the circle as follows.
First move
$\s_j$ to the left of $\Delta_n$, when it becomes $\s_{n-j}$. Then move it round the circle to the
end of the word. It can then be moved left past all the remaining generators of $\alpha^{-1}$, and
past
$\Delta_n$ once more, to become $\s_j$. Finally move this round the circle to the right-hand end of
$\alpha^{-1}$.

Consequently $\be$ is conjugate to $\Delta_n \sigma_1^{-1}\sigma_2^{-1}\cdots\sigma_k^{-1}$.
\end{proof}

\section{Non-conjugacy results}
When $n=6$ it is possible to have two positive permutation braids with the same number of crossings
which close to different knots. The permutations $(124536)$, with braid
$\s_3\s_4\s_3\s_2\s_5\s_4\s_3\s_2\s_1$, and $(132546)$, with braid 
$\s_2\s_1\s_4\s_3\s_5\s_4\s_3\s_2\s_1$, close to the $(2,5)$ torus knot and the sum of two trefoils
respectively, so writhe no longer determines conjugacy.

In \cite{Hadji}, Hadji gave examples of two non-conjugate
positive permutation braids with $n=16$ each closing to the same connected sum of three knots.

In fact, non-conjugate positive permutation braids which close to the same knot show up first when
$n=6$.
\begin{theorem} \label{nonconj}
The positive permutation braids in $S_6^+$ with permutations $(165324)$ and $(152643)$ have the same
closure but are not conjugate.
\end{theorem}
\begin{proof}
The braids, shown below,  can be written $\be=\sigma_1\s_3\s_5\s_2\s_4\s_1\s_3\s_2\s_1$ and
$\gamma=\s_2\s_4\s_3\s_5\s_2\s_4\s_1\s_3\s_2$ respectively. Both of these can be reduced by Markov
moves to the $4$-braid $\s_1\s_3\s_2\s_1\s_3\s_2\s_1$, so both close to the $(2,5)$ torus knot.

$$\be\quad=\quad\Firstbraid\qquad \gamma\quad=\quad \Secondbraid$$

The squares of the two braids are shown here with the strings which form one component of the
closure emphasised.

$$\be^2\quad=\quad\Firstsquare\qquad \gamma^2\quad=\quad \Secondsquare$$

If $\be$ and $\gamma$ are  conjugate then so are $\be^2$ and $\gamma^2$.
 Now the closure of $\be^2$ is a link with two components, each of which turns
out to be the trefoil knot, while the two components of the closure of $\gamma^2$ are trivial knots.
Hence $\be^2$ and $\gamma^2$ are not conjugate.
\end{proof}

An alternative check can be made by calculating the 2-variable Alexander polynomial of the  
link consisting of the closure of $\be$ and its axis. If $\be$ is conjugate to $\gamma$ this link
is isotopic to the closure of $\gamma$ and its axis.  Its polynomial is in general the
characteristic polynomial of the reduced Burau matrix of the braid,
\cite{Morton}.  For $\be$ above, the polynomial is $$t^9x^5+t^7x^4+t^5x^3+t^4x^2+t^2x+1,$$
which differs, up to
 multiples of $\pm t^ix^j$, from the polynomial 
$$t^9x^5+t^7x^4+(2t^5-t^4)x^3+(2t^4-x^5)x^2+t^2x+1$$
for $\gamma$. 

\bigskip
Other tests for conjugacy, which also rely in effect on invariants of a closed braid in a solid
torus, can be used to give a contrasting result to theorem \ref{2} about positive braids which
close to the trefoil, when we do not restrict to positive permutation braids.

\begin{theorem}\label{trefoil} If $\be\in B_n^+$ closes to the trefoil knot then $\be$ is conjugate
to
$\be(i)=\s_1\s_2\ldots\s_{i-1}\s_i^3\s_{i+1}\ldots\s_{n-1}$ for some $i$. Two such braids $\be(i),
\be(k)$ are conjugate if and only if $k=i$ or $k=n-i$.
\end{theorem}
\begin{remark} When $n=4$ the braids are examples of the construction of Murasugi and Thomas,
\cite{Thomas}. They show that the braids $\s_1^p\s_2^q\s_3^r$ and  $\s_1^p\s_2^r\s_3^q$, with
$p,q,r$ odd, which close to isotopic knots, are not conjugate when $q\ne r$. Their proof uses the
exceptional homomorphism from $B_4$ to $B_3$ defined by $\s_1,\s_3\mapsto\s_1,\s_2\mapsto\s_2$,
observing that the braids map to $\s_1^{p+r}\s_2^q$ and $\s_1^{p+q}\s_2^r$, which close to links
with different linking numbers.
\end{remark}

\begin{proof}[of theorem \ref{trefoil}]

\noindent
1. {\it Conjugacy}. \quad The only difference from the argument of theorem \ref{2} is that one
generator
$\s_i$ may occur three times, with $\s_{i-1}$ and $\s_{i+1}$ both lying on the circle between the
same pair of occurrences of $\s_i$. Then all three  occurrences of $\s_i$ can be moved together
and remain as a block on the circle, while the other generators are put in consecutive order, as in
theorem \ref{1}. This shows that every such braid is conjugate to some $\be(i)$.
To see that the braids $\be(i)$ and $\be(n-i)$ are conjugate, first conjugate $\be(i)$ by
$\Delta_n$, taking $\s_i^3$ to $\s_{n-i}^3$ and then rearrange as above.

\noindent
2. {\it Non-conjugacy}. \quad Any closed braid represents an element in the framed Homfly skein of
closed braids in the annulus \cite{Morton2}. The closure of $\s_1\s_2\ldots\s_{k-1}$ represents an
element
$A_k$. The skein itself admits a commutative product, represented by the closures of split braids.
The subspace spanned by the closure of braids in $B_n$ has a basis consisting of monomials
$A_{i_1}\ldots A_{i_k}$ with $i_1+\cdots+i_k=n$. Coefficients in the skein can be taken as integer
polynomials in a variable $z$. In the Homfly skein of braids before closure, we have
$\s_i^3=c(z)\s_i+d(z)$, for some fixed non-zero polynomials $c(z),d(z)$, so that
$\be(i)=c(z)\s_1\ldots\s_{n-1}+d(z)\s_1\ldots\s_{i-1}\s_{i+1}\ldots\s_{n-1}$ in this skein. Its
closure then represents $c(z)A_n+d(z)A_iA_{n-i}$ in the skein of the annulus.

If $\be(i)$ and $\be(k)$ are conjugate then they have the same closure in the annulus. Then
$$c(z)A_n+d(z)A_iA_{n-i}=c(z)A_n+d(z)A_kA_{n-k},$$ and hence $A_iA_{n-i}=A_kA_{n-k}$. The
monomials form a basis in the skein of the annulus, so $k=i$ or $k=n-i$.

\end{proof}

\begin{remark}
This same calculation can be used to show that the Conway polynomial of the closure of $\be(i)$ and
its axis differs from that of $\be(k)$ and its axis except when $k=i$ or $n-i$.
\end{remark}

\section{Conjugacy classes for $6$ and more strings}

We have a short Maple procedure to list the positive permutation braids on $n$ strings which close
to knots, according to their number of crossings. 
\medskip

\noindent
{\it The case $n=6$.}

When $n=6$ this list contains $16$ positive permutation braids with $5$
crossings, $32$ with $7$ crossings, $44$ with $9$ crossings, $22$ with $11$ crossings and $6$ with
$13$ crossings.

By theorems \ref{1}, \ref{2} and \ref{large}  those with $5,7$ or
$13$ crossings form complete conjugacy classes, 
and represent the trivial knot, the trefoil and the $(3,5)$ torus knot respectively. An inductive
count shows that there are in general $2^{n-2}$ braids in $S_n^+$ which represent the trivial knot.

Among the braids with $9$ crossings there is one conjugacy class consisting of $4$ braids which
close to the sum of two trefoils, and  two classes of braids which close to the $(2,5)$ torus knot.
There are just $2$ braids, $\gamma=\s_2\s_4\s_3\s_5\s_2\s_4\s_1\s_3\s_2$ and its conjugate by the
half-twist, in the conjugacy class of the braid
$\gamma$ discussed in theorem \ref{nonconj}, while the remaining $38$ braids are
conjugate to
$\be=\sigma_1\s_3\s_5\s_2\s_4\s_1\s_3\s_2\s_1$.

The braids with $11$ crossings fall into two conjugacy classes, one containing $6$ braids which
close to the $(3,4)$ torus knot, and the other containing  $16$ braids which close to the $(2,7)$
torus knot.

\medskip

\noindent
{\it The case $n=7$.}

 When $n=7$ there are $32$ positive permutation braids with $6$ crossings,
$88$ with $8$ crossings,
$176$ with
$10$ crossings,
$202$ with $12$ crossings, $134$ with $14$ crossings, $70$ with $16$ crossings and $18$ with $18$
crossings.

 Again those with $6,8$ or $18$ crossings represent complete conjugacy classes; we have
not attempted to analyse the other classes any further, or to consider in detail any cases where
$n>7$.

\noindent
Original version September 2003. Current version December 2003.
\end{document}